\documentclass[12pt,draft] {iopart}

\usepackage[english]{babel}
\usepackage{graphicx,color}
\usepackage{amsfonts}
\usepackage{amssymb,amsthm,scalefnt,mathrsfs}
\usepackage[all]{xy}

\newtheorem{theo}{Theorem}[section]
\newtheorem{cor}[theo]{Corollary}

\newtheorem{defn}[theo]{Definition}
\newtheorem{rem}[theo]{Remark}
\newtheorem{ex}[theo]{Example}

\newcommand{\abs}[1]{\left\vert#1\right\vert}
\newcommand{\set}[1]{\left\{#1\right\}}
\newcommand{\C}[1]{\left( #1\right)}
\newcommand{\CC}[1]{\left[#1\right]}

\newcommand{\A}{\mathbf{A}}
\newcommand{\B}{\mathbf{B}}

\newcommand{\X}{\mathbf{X}}

\newcommand{\Y}{\mathbf{Y}}
\newcommand{\Prob}{\mathbb{P}}
\newcommand{\E}{\mathbb{E}}
\newcommand{\R}{\mathbb{R}}

\newcommand{\Z}{\mathbb{Z}}
\newcommand{\N}{\mathbb{N}}

\newcommand{\cqd}{\begin{flushright} $\square$\\\end{flushright}}

\begin{document}

\title[Dualities for Probabilistic Cellular Automata]{Dualities for Multi-State Probabilistic Cellular Automata}

\author{F. Javier L\'{o}pez$^1$,
Gerardo Sanz$^1$ and Marcelo Sobottka$^2$}

\address{$^1$Departamento de M\'{e}todos Estad\'{\i}sticos and BIFI, Facultad de
Ciencias, Universidad de Zaragoza, C/ Pedro Cerbuna 12, 50009
Zaragoza, Spain.}

\address{$^2$Departamento de Ingenier\'{\i}a Matem\'atica, Facultad de
Ciencias F\'{\i}sicas y Matem\'{a}ticas, Universidad de Concepci\'{o}n, Casilla
160-C, Correo 3, Concepci\'{o}n, Chile.}

\eads{\mailto{javier.lopez@unizar.es},
\mailto{gerardo.sanz@unizar.es}, \mailto{msobottka@udec.cl}}

\date{\ }



   \setcounter{secnumdepth}{2}
   \setcounter{tocdepth}{1}   


\begin{abstract}

In this paper a new form of duality for probabilistic cellular
automata (PCA) is introduced. Using this duality, an ergodicity
result for processes having a dual is proved. Also, conditions on
the probabilities defining the evolution of the processes for the
existence of a dual are provided. The results are applied to wide
classes of PCA which include multi-opinion voter models, competition
models and the Domany-Kinzel model.

\end{abstract}


\bigskip
\hrule
\noindent
{\footnotesize\em This is a pre-copy-editing, author-produced PDF of an article accepted for publication in Journal of Statistical Mechanics, following peer review. The definitive publisher-authenticated version {\em F. J. L\'{o}pez, G. Sanz and M. Sobottka. Dualities for multi-state probabilistic cellular automata. Journal of Statistical Mechanics (2008),  P05006, doi:10.1088/1742-5468/2008/05/P05006}, is available online at:\break http://www.iop.org/EJ/abstract/1742-5468/2008/05/P05006 .}
\hrule
\bigskip

\section{Introduction}

Probabilistic Cellular Automata (PCA) are discrete-time stochastic
processes with state space $\X:=W^{\Z^d}$, where
$W=\set{1,\ldots,M}$ and $d\geq 1$, with finite range interactions
on the integer lattice $\Z^d$ (this is interpreted as at each point
$z\in Z^d$ there is a particle which  take values in the set $W$;
the points $x\in{\bf X}$ are called configurations). More precisely,
suppose $N\subset\Z^d$ is a finite subset and $f:W^N\times W\to
[0,1]$ is a transition function, then a PCA is a discrete-time
homogenous Markov process $\eta_t=\{\eta_t(z)\in W:\ z\in\Z^d\}$
whose evolution is given, for all $s\in\N$, $z\in\Z^d$, $w\in W$,
and $x\in\X$, by
\begin{equation}\label{PCA}\Prob\set{\eta_{s+1}(z)=w|\
\eta_s=x}=f\bigl((x(z+n))_{n\in N},w\bigr).\end{equation} In other
words, $\eta_t$ is an interacting particle system on $\X$ in
discrete time such that the transitions at different sites are
conditionally independent, given the current state of the system.
For definitions and main results on PCA, the reader can consult
\cite{leb} and \cite{toom}.

The set $\X=W^{\Z^d}$ is endowed with the product topology. Let
  ${\cal P}(\X)$ be the set of probability measures on
$\X$ with the weak* convergence topology. Note that
  ${\cal P}(\X)$ is compact with respect to this
topology.

Let $\eta_t$ and $\xi_t$ be two stochastic processes (which can
evolve either in continuous or discrete time) with state spaces $\X$
and $\Y$ respectively, and suppose $H:\X\times \Y\to\R$ is a bounded
measurable function, then $\eta_t$ and $\xi_t$ are dual to one
another with respect to $H$ if for all $x\in \X$ and $y\in \Y$
(Definition II.3.1 of \cite{liggett1})
\begin{equation}\label{dual}\E_{\eta_0=x}\left[H(\eta_s,y)\right]
=\E_{\xi_0=y}\left[H(x,\xi_s)\right],\end{equation} where
$\E_{\gamma_0=\gamma}[\cdot]$ denotes the expectation when the
process $\gamma_t$ starts with configuration $\gamma$.

If $\eta_t$ and $\xi_t$ are both discrete-time Markov chains with
transition matrices $P$ and $Q$ respectively, then the previous
expression for duality can be written as
\begin{equation}\label{markov}P^sH=H(Q^s)^T,\end{equation}
where the superscript $^T$ stands for the transpose. Furthermore,
due to the Markovian property, if the duality equation holds for
$s=1$ then it holds for all $s\in\N$.

Duality allows obtaining relevant information about the evolution of
$\eta_t$, which usually has uncountable state space, by studying the
evolution of $\xi_t$ which in general is chosen having countable
state space. In fact, duality has been widely used in the study of
spin systems (continuous-time interacting particle systems with
$W=\{0,1\}$), where $\Y$ is the set of all finite subsets of $\Z^d$
(see \cite{liggett1}, \cite{Dur}, \cite{gray}, and \cite{katori1}).
For the case where $\eta_t$ is a continuous-time interacting
particle system with $\abs{W}\geq 2$, L\'{o}pez and Sanz \cite{ls}
studied an extended version of (\ref{dual}), which is given by the
equation:
\begin{equation}\label{dualV}\E_{\eta_0=x}\left[H(\eta_s,y)\right]=\E_{\xi_0=y}
\left[H(x,\xi_s) e^{-\int_0^s{V(\xi_u)du}}\right],\end{equation} for
all $x\in \X$ and $y\in \Y$, and for some function $V:\Y\to
[0,\infty)$. In particular, this duality equation was used in
\cite{ls} to obtain results about the long time behaviour of the
continuous-time multi-opinion noisy voter model and the
continuous-time 3-opinion
noisy biased voter model.\\

For the discrete-time case, Katori et al. \cite{katori2} develop a
theory of duality for two state PCA based on the equation
\begin{equation}\label{eqDK}\E\left[\mathfrak{c}^{\vert\xi_s^B\cap
A\vert}\right]=\E\left[\mathfrak{c}^{\vert\eta_s^A\cap
B\vert}\right],\end{equation} where $\mathfrak{c}$ is a parameter
depending on the processes $\eta_t$ and $\xi_t$; $A,B\subseteq\Z$,
with at least one of them being finite; $\eta_s^A$ and $\xi_s^B$
represent the set of particles with value 1 at time $s$ starting
from configurations where only the particles in $A$ and $B$,
respectively, have the value 1; and $\vert\cdot\vert$ stands for the
cardinality of a set (this equation was first considered for
interacting particle systems by Sudbury and Lloyd in \cite{SL}). In
\cite{katori2}, this relationship is used to study duality for some
cases of the Domany-Kinzel model, which is the PCA with state space
$\{0,1\}^\mathbb{Z}$, and evolution defined by the parameters
$a_0,a_1,a_2\in [0,1]$ as follows:
\begin{equation}\label{DK}\hspace{-1cm}\Prob\set{\eta_{s+1}(z)=1|\
\eta_s(z-1)=w_{-1},\eta_s(z)=w_{0},
\eta_s(z+1)=w_{1}}=a_{w_{-1}+w_1}\end{equation}

When $a_0=0$ and $a_2\leq a_1$, Katori et al. \cite{katori2} have
presented a dual for $\eta_t$. Such dual is a model in which each
stage is updated by a $a_1$-thinning of all sites, followed by an
application of the Domany-Kinzel model with parameters $a_0'=0$,
$a_1'=1$, and $a_2'=a_2/a_1$. When $a_1>a_2$, they also
  studied the limit behaviour of the process. The
technique used in \cite{katori2} was to locate the process $\eta_t$
in a finite subset of $\Z$ to represent it by a transition matrix
$P$, and so to find $H$ and $Q$ which satisfy the duality equation
(\ref{markov}). The technique used by Katori et al. depends strongly
on the fact that $a_0=0$ (see also chapter 5 of \cite{Dur} for
results on this model, with $a_0=0$, using percolation theory). In
\cite{konno1} and \cite{konnoma}, those results are extended to a
larger class of one dimensional two state PCA.

For multi-state PCA, an equation analogous to (\ref{eqDK}) is
introduced in \cite{konno2} to define self-dual processes (that is,
processes $\eta_t$ that satisfy the equation with $\xi_t=\eta_t$)
and give conditions for self-duality.

Our objective in this paper is to give a general theory of duality
for multi-state PCA. The duality equation we define (see
(\ref{duald}) below) is given for general duality functions $H$ and
is not restricted to self-duality. Moreover, the introduction of an
auxiliary function $d$ in the equation allows the inclusion of wider
classes of processes for each dual function $H$. We will show the
usefulness of our concept by obtaining an ergodicity result (Theorem
\ref{theo(H,d)}) for processes satisfying (\ref{duald}). Equation
(\ref{duald}) and Theorem \ref{theo(H,d)} are stated for general $H$
and we show their applicability with two particular cases of $H$
which cover two wide classes of PCA, providing conditions on their
transition probabilities such that they satisfy (\ref{duald}) and,
in that case, providing conditions for their ergodicity (that is,
the existence of an unique equilibrium measure and the convergence
of the process to that measure). We will also apply our results to
some examples.

The duality equation we propose is

{\defn Given two discrete-time Markov processes, $\eta_t$ with state
space $\X$ and $\xi_t$ with state space $\Y$, and $H:\X\times
\Y\to\mathbb{R}$ and $d:\Y\to [0,\infty)$ bounded measurable
functions, we say $\eta_t$ and $\xi_t$ are dual to one another with
respect to $(H,d)$ if
\begin{equation}\label{duald}\E_{\eta_0=x}\left[H(\eta_1,y)\right]
=d(y)\E_{\xi_0=y}\left[H(x,\xi_1) \right].\end{equation}}

Notice that the duality equation with respect to $(H,d)$ is only
written for the one-step evolution of the processes $\eta_t$ and
$\xi_t$ as they are both Markov processes, but an induction argument
yields the duality equation for the $s$ step evolution of the
process:

\begin{equation}\label{duald_n}E_{\eta_0=x}[H(\eta_s,y)]=
d(y)E_{\xi_0=y}[d(\xi_1)\cdots
d(\xi_{s-1})H(x,\xi_s)].\end{equation}

When $\eta_t$ and $\xi_t$ are Markov chains with transition matrix
$P$ and $Q$, respectively, equation (\ref{duald}) can be written as
\begin{equation}\label{dualD}PH=H(DQ)^T,\end{equation} where $D$
 is the diagonal matrix
with $D_{yy}=d(y)$. In this case, the general expression for time
$s\geq 1$ is $P^sH=H\left((DQ)^s\right)^T$.

The advantage of considering the function $d$ on the right side of
(\ref{duald}) is that processes that do not have dual with respect
to some function $H$ in (\ref{dual}) may have it with respect to
$(H,d)$ (see Remarks \ref{d_relevance1} and \ref{d_relevance2}
below). When $d\equiv 1$, we have the classical notion of duality
(\ref{dual}), which means that the evolution of $\eta_t$ is easily
understood from the evolution of $\xi_t$. When $d\not\equiv 1$, the
relationship between $\eta_t$ and $\xi_t$ becomes more complicated;
nevertheless, in that case, important information about $\eta_t$ can
be obtained from $\xi_t$.

Although equations (\ref{duald}) and (\ref{duald_n}) can be seen as
a discrete counterpart of equation (\ref{dualV}), given for the
continuous-time case, there are several differences between the
techniques used to exploit the former and the latter. In fact,
instead of using directly equation (\ref{duald_n}) to prove results
on ergodicity for the process, we show that there is a
transformation $\tilde \xi_t$ of the dual process $\xi_t$ such that
$\eta_t$ and $\tilde\xi_t$ are dual in the classical sense (with
respect to a modified function $\tilde H$), see Theorem
\ref{theo(H,d)}(a) below. This approach is useful, for instance, for
writing down the explicit form of the invariant measure
(\ref{forminv}) and for providing alternative conditions for
ergodicity (see Remark \ref{remtau}) which are used, for instance,
in the study of the Domany-Kinzel model (Corollary \ref{cordkm}).

The organization of the paper is as follows. In {\S}\ref{general} we
present the basic tools for the duality of stochastic processes and
obtain conditions on $H$ and $d$ for the ergodicity of PCA. In
{\S}\ref{examples} we analyze, using the duality theory developed, two
general classes of PCA. For each of these classes, we first obtain
conditions on the transition probabilities of
 the processes so that they have a dual and then give
conditions for their egodicity. The results are applied to specific
models such as the multi-opinion noisy voter model, the
Domany-Kinzel model and competition models. For the particular case
of the Domany-Kinzel model with parameters $0\leq a_0\leq a_1\leq
a_2$ we find a dual process for many cases not covered by \cite{Dur}
or \cite{katori2}; and we get new results for the ergodicity of this
model.

\section{Dualities}\label{general}

In this section we  consider probabilistic cellular automata
$\eta_t$ for which there exist a pair of functions $(H,d)$ and a
dual process $\xi_t$ with respect to $(H,d)$. Sufficient conditions
on $(H,d)$ for the ergodicity shall be presented together with the
characterization of the corresponding equilibrium
measure.\\

In what follows, we first recall some basic definitions and results
about dualities for stochastic processes. Suppose $\eta_t$ on $\X$
and $\xi_t$ on $\Y$ are two stochastic processes that are dual to
one another with respect to a function $H$. Given a probability
measure $\mu$ on $\X$, define for any $y\in\Y$

\begin{equation*}\hat{\mu}(y):=\int_{\X}H(x,y)d\mu(x).\end{equation*}

For any $s\geq 0$, denote by $\mu_s$ the distribution obtained from
the initial distribution $\mu$ when $\eta_t$ evolves until time $s$.
In our setting, $\X=W^{\Z^d}$ and $\eta_t$ is a PCA given by
(\ref{PCA}); then $\mu_s$ is defined on the cylinder
$\mathbf{w}=[w_{z_1},\ldots,w_{z_m}]:=\set{x\in\X:\ x(z_i)=w_{z_i},
1\leq i\leq m}$ by
$$\mu_s(\mathbf{w})=\sum_{\mathbf{u}\in\mathcal{C}_s}
\mu(\mathbf{u})\Prob_{\eta_0\in\mathbf{u}}\set{\eta_s\in\mathbf{w}},$$
where $\mathcal{C}_s$ is the family of all cylinders of $\X$
defined on the coordinates $\set{z_i+n_j:\ i=1,\ldots,m;\ n_j\in
N}$\sloppy.

Therefore, we have that
$$\hat{\mu}_s(y)=\int_{\X}H(x,y)d\mu_s(x)=
\int_{\X}\E_{\eta_0=x}\CC{H(\eta_s,y)}d\mu(x),$$ and by the duality
equation (\ref{dual}) it follows that
\begin{equation}\label{mu_s}\begin{array}{lcl}\hat{\mu}_s(y) & = &
\int_{\X}\E_{\eta_0=x}\CC{H(\eta_s,y)}d\mu(x)\\\\
& = & \int_{\X}\E_{\xi_0=y}\CC{H(x,\xi_s)}d\mu(x)\\\\
& = & \E_{\xi_0=y}\CC{\int_{\X}H(x,\xi_s)d\mu(x)}\\\\
& = & \E_{\xi_0=y}\CC{ \hat{\mu}(\xi_s)}.
\end{array}\end{equation}

Suppose the set of the linear combinations of the functions
$\set{H(\cdot,y):\ y\in\Y}$ is dense in $C(\X)$, the space of
continuous real functions on $\X$. If for any $y$,
$E_{\xi_0=y}[\hat\mu(\xi_s)]$ converges to some function of $y$ (not
depending on $\mu$) as $s$ goes to $\infty$, then for any $f\in
C(\X)$, it follows that $\int_\X f d\mu_s$ converges to some
functional $\tilde{\nu}(f)$, which means $\mu_s$ converges in the
weak* topology to some probability measure $\nu$. Furthermore, if
the above holds for any initial distribution $\mu$, it implies the
ergodicity of $\eta_t$: that is, the convergence in distribution of
the process to the equilibrium measure $\nu$ from any initial
measure.

Now we are able to study the case of the duality with respect to
$(H,d)$.\\

{\theo\label{theo(H,d)} Suppose $\eta_t$ is a Markov process with
state space $\X$ and $\xi_t$ is a Markov chain with state space
$\Y$, which are dual to one another with respect to $(H,d)$, where
$0\le d(y)\leq 1$ for all $y\in\Y$. Then:

\begin{description}

\item[(a)] There exist a Markov chain $\tilde{\xi}_t$ with state space
$\tilde{\Y}=\Y\cup\{\wp\}$, where $\wp$ is a new state, and a
bounded measurable function $\tilde{H}:\X\times \tilde{\Y}\to\R$
such that $\eta_t$ and $\tilde{\xi}_t$ are dual to one another with
respect to $\tilde{H}$;

\item[(b)] denoting by $\Theta$ the set of all absorbing states of
$\xi_t$, if

\begin{enumerate}

\item the linear combinations of $\bigl\{H(\cdot,y):\ y\in\Y\bigr\}$ is a dense set of
$C(\X)$;

\item $d(y)< 1$ for any $y\notin\Theta$,
and $\displaystyle \sup_{y\in\Y:\ d(y)<1}\ d(y)<1$;

\item $H(\cdot,\theta)\equiv c(\theta)$ for all
$\theta\in\Theta$ with $d(\theta)=1$;

\end{enumerate}

then $\eta_t$ is ergodic and its unique equilibrium measure is
determined for any $y\in\Y$ by
\begin{equation}\label{forminv}
\hat{\nu}(y)=\sum_{\begin{array}{c}\theta\in\Theta,\\d(\theta)=1\end{array}}c(\theta)\Prob_{\tilde{\xi}_0=y}\set{\tilde{\xi}_\tau=\theta},
\end{equation}
where $\tau$ is the hitting time of $\{\theta\in\Theta:\
d(\theta)=1\}\cup\{\wp\}$ for $\tilde{\xi}_t$.

\end{description}
}

\proof {\bf (a)} Since $0\leq d(y)\leq 1$ for all $y\in\Y$, we can
define a Markov chain $\tilde{\xi}_t$ with state space
$\tilde{\Y}=\Y\cup\set{\wp}$ and transition probabilities given by
\begin{equation}\label{tilde xi_t}
\Prob_{\tilde{\xi}_0=\tilde{y}}\{\tilde{\xi}_1=\tilde{y}'\}:=\left\{\begin{array}{lcl}
d({\tilde{y}}) \Prob_{{\xi}_0=\tilde{y}}\{{\xi}_1=\tilde{y}'\} & {,\ if} & \tilde{y},\tilde{y}'\in \Y,\\\\
1-d(\tilde{y}) & {,\ if} & \tilde{y}\in \Y,\ \tilde{y}'=\wp\\\\
1  & {,\ if} & \tilde{y}=\tilde{y}'=\wp\end{array}\right.
\end{equation}

Note that the new state $\wp$ is an absorbing state of the process
$\tilde{\xi}_t$.

Define $\tilde{H}:\X\times\tilde{\Y}\to\R$ by
\begin{equation*}
\tilde{H}(x,\tilde{y}):=\left\{\begin{array}{lcl}\displaystyle
H(x,\tilde{y}) & {,\ if} & \tilde{y}\in \Y,\\\\
0 & {,\ if} & \tilde{y}=\wp \end{array}\right.
\end{equation*}

It is trivial that (\ref{dual}) holds for $x\in\X$ and $\wp$ since
$\tilde{H}(x,\wp)=0$. On the other hand, if $\tilde{y}\in\Y$, then
$$\begin{array}{lll}
\E_{\eta_0=x}\left[\tilde{H}(\eta_1,\tilde{y})\right] & = &
\displaystyle \E_{\eta_0=x}\left[H(\eta_1,\tilde{y})\right] \\\\ &=&
d(\tilde{y})\E_{\xi_0=\tilde{y}}\left[H(x,\xi_1)\right]\\\\
& = & \displaystyle d(\tilde{y})\sum_{y'\in \Y}H(x,y')\Prob_{\xi_0=\tilde{y}}\set{\xi_1=y'}\\\\
& = &  \displaystyle \sum_{\tilde{y}'\in\Y\cup\set{\wp}}\tilde
H(x,\tilde{y}')\Prob_{\tilde\xi_0=\tilde{y}}\set{\tilde\xi_1=\tilde{y}'}\\\\
& = & \displaystyle
\E_{\tilde{\xi}_0=\tilde{y}}\left[\tilde{H}(x,\tilde{\xi}_1)\right].
\end{array}$$

Therefore, $\eta_t$ and $\tilde{\xi}_t$ are dual to one another with
respect to $\tilde{H}$.

{\bf (b)} Note first that by {\em i}, the set of linear combinations
of $\tilde H(\cdot,\tilde{y})$, with $\tilde{y}\in\tilde\Y$, is
dense on $C(\X)$.

Denote by $\tilde{\Theta}:=\set{\theta\in\Theta:\
d(\theta)=1}\cup\set{\wp}$ the set of all absorbing states of
$\tilde{\xi}_t$ and denote by $\tau$ the hitting time of
$\tilde{\Theta}$ for $\tilde{\xi}_t$, that is,
$$\tau:=\min\set{s\geq 0:\quad \tilde{\xi}_s\in \tilde{\Theta}}.$$

From hypothesis {\em ii}, we obtain that there exists $a<1$ such
that for any non-absorbing state $\tilde{y}$ of $\tilde{\xi}_t$ we
have that $d(\tilde{y})\leq a$ and
$\Prob_{\tilde\xi_0=\tilde{y}}\set{\tilde\xi_1\not\in\tilde\Theta}
\leq a$. It follows that for any
$\tilde{y}\in\tilde{\Y}\setminus\tilde{\Theta}$, and any $s\geq 1$:
\begin{equation*}
\begin{array}{lcl}
\Prob_{\tilde{\xi}_0=\tilde{y}}\set{\tilde\xi_{s}\not\in\tilde\Theta}
&=&\displaystyle\sum_{\tilde{z}\in\tilde{\Y}\setminus\tilde{\Theta}}
\Prob_{\tilde{\xi}_0=\tilde{y}}\set{\tilde\xi_{s-1}=\tilde{z}}
\Prob_{\tilde\xi_{s-1}=\tilde{z}}\set{\tilde\xi_{s}\not\in\tilde\Theta}
\\\\
& \leq &\displaystyle a
\sum_{\tilde{z}\in\tilde{\Y}\setminus\tilde{\Theta}}
\Prob_{\tilde{\xi}_0=\tilde{y}}\set{\tilde\xi_{s-1}=\tilde{z}}
\\\\
&=&
a\Prob_{\tilde\xi_0=\tilde{y}}\set{\tilde\xi_{s-1}\not\in\tilde\Theta}.
\end{array}
\end{equation*}

Therefore, by induction and noting that
 $\Prob_{\tilde{\xi}_0=\tilde{y}}\set{\tau=\infty}  \le
\Prob_{\tilde{\xi}_0=\tilde{y}}\set{\tilde\xi_{s}\not\in\tilde\Theta}$
for any $s\ge1$, we get
$\Prob_{\tilde{\xi}_0=\tilde{y}}\set{\tau=\infty} \leq a^s$ and,
since it holds for all $s\ge 1$, we deduce
$\Prob_{\tilde{\xi}_0=\tilde{y}}\set{\tau<\infty}=1$.\\

Defining, for any probability measure $\mu\in{\cal P}(\X)$ and any
$\tilde y\in\tilde\Y$, $ \hat\mu(\tilde y):=\int_{\X}\tilde
H(x,\tilde y)d\mu(x)$, it follows that for any $\tilde y\in\tilde\Y$
 \begin{equation*}\hspace{-1cm}
\begin{array}{lcl}\displaystyle\lim_{s\rightarrow\infty}\hat{\mu}_s(\tilde y) & =
& \displaystyle \lim_{s\rightarrow\infty}\E_{\tilde{\xi}_0=\tilde{y}}\CC{ \hat{\mu}(\tilde{\xi}_s)}\\\\
& = &\displaystyle\lim_{s\rightarrow\infty}
\sum_{\tilde{\theta}\in\tilde{\Theta}}\E_{\tilde{\xi}_0=\tilde
y}\CC{ \hat{\mu}(\tilde{\xi}_s)| \ \tilde{\xi}_\tau=\tilde{\theta},\
\tau\leq s }\Prob_{\tilde{\xi}_0={\tilde
y}}\set{\tilde{\xi}_\tau=\tilde{\theta},\ \tau\leq
s}\\\\
&& \displaystyle
+\lim_{s\rightarrow\infty}\E_{\tilde{\xi}_0=\tilde{y}}\CC{
\hat{\mu}(\tilde{\xi}_s)| \ \tau>
s}\Prob_{\tilde{\xi}_0=\tilde{y}}\set{\tau>
s}\\\\
& = &\displaystyle \sum_{\theta\in\Theta,\
d(\theta)=1}c(\theta)\Prob_{\tilde{\xi}_0=\tilde
y}\set{\tilde{\xi}_\tau=\theta},
\end{array}\end{equation*}
where the first equality follows from the fact that $\eta_t$ and
$\tilde{\xi}_t$ are dual to one another with respect to $\tilde{H}$
and equation (\ref{mu_s}), and the last one is due to {\em iii},
$\tilde{H}(\cdot,\wp)\equiv 0$ and
$\tilde{\Theta}=\set{\theta\in\Theta:\ d(\theta)=1}\cup\set{\wp}$.
\cqd

{\rem \label{remtau} From the proof of Theorem \ref{theo(H,d)}, it
is straightforward to check that condition {\em ii} can be replaced
by $\Prob_{\xi_0=y}\set{\tau'<\infty}=1$ for all $y\in\Y$, where
$\tau'$ is the hitting
time of $\Theta$ for $\xi_t$, since this implies $\Prob_{\tilde\xi_0=\tilde y}\set{\tau<\infty}=1$ for all $\tilde y\in\tilde Y$.}\\

{\rem As equation (\ref{duald}) can be seen as a discrete-time
counterpart of (\ref{dualV}), it is interesting to compare the form
of the limit measure (\ref{forminv}) with the expression given in
(2.4) of \cite{ls} for continuous time:
$$\hat\nu(y)=\E_{\xi_0=y}\left[e^{-\int_0^\tau V(\xi_s)ds}{\mathbf
1}_{\{\tau<\infty\}}\right].$$ Note that, in \cite{ls}, the
expression is given through the dual process $\xi_t$, studying its
complete evolution from zero time to the time of absorption, while
here it is given only through the absorbing point of the modified
process
$\tilde\xi_t$ instead of the original dual process $\xi_t$.}\\

{\rem Theorem \ref{theo(H,d)} gives conditions for the ergodicity of
the process and, in that case, the form of the limit measure. Note
that, even when the process is not ergodic, equation (\ref{duald})
and the computations made in the proof of the theorem can be used to
study the behaviour of $\eta_t$ as a function of the evolution of
$\tilde\xi_t$.}

\section{Nearest-neighbour probabilistic cellular automata}\label{examples}

In this section we will use the duality equation (\ref{duald}) to
study conditions for the ergodicity and the limit behaviour of some
PCA models evolving on $\Z$.

First, let us introduce a definition:

{\defn We say that a PCA $\eta_t$ with state space
$\X=W^\Z=\{1,\ldots,M\}^{\Z}$ is a nearest-neighbour PCA if
$N=\{-1,0,1\}$.}

Thus, the evolution of such a process is given for any $i,j,k,m\in
W$, $z\in\Z$ and $s\in\N$ by
$$\Prob\set{\eta_{s+1}(z)=m|\ \eta_s(z-1)=i,\eta_s(z)=j,
\eta_s(z+1)=k}=p_{ijk,m}.$$

We will treat two wide classes of nearest-neighbour PCA, including
classical voter models, linear and non-linear voter models, noisy
voter models, biased voter models and some competition models.

\subsection{Multi-opinion voter models}\label{ssmonvm}

  Throughout this subsection, $\eta_t$ will be a
nearest-neighbour PCA with $M$ states, and the state space of the
dual processes will be $\Y= Y^{M-1}$ (where $Y$ is the collection
of finite subsets of $\Z$).

Let $H:\X\times\Y\to\R$ be the function defined for all $x=\{x(z)\in
W\}_{z\in\Z}\in\X$ and $\A=(A_i)_{1\leq i\leq M-1}\in\Y$
 by
\begin{equation}\label{H1} H(x,\A)=\left\{\begin{array}{lll} 1  &,&\ {if}\
x(z)=i,\ \forall z\in A_i,\ \forall i=1,\ldots, M-1\\
0 &,&\ {otherwise}\end{array}\right..\end{equation}

Our aim is to construct a Markov chain $\A_t$ with state space $\Y$
and a function $d:\Y\to [0,\infty)$ that verify equation
(\ref{duald}) with $H$ defined in (\ref{H1}).

Suppose $\eta_t$ and $\A_t=(A_{1,t},\ldots, A_{M-1,t})$ are dual to
one another with respect to $(H,d)$, for some function $d$, and let
$x\in\X$ and $\A=(A_i)_{1\leq i \leq M-1}\in\Y$, where
$\set{A_i}_{1\leq i \leq M-1}$ is a pairwise disjoint family of
subsets of $\Z$. Then,  equation (\ref{duald}) can be written as:

\footnotesize
\begin{equation}\label{MONVM}\hspace{-2.5cm}\begin{array}{rl} \Prob_{\eta_0=x}\set{\eta_1(z)=i,\forall
z\in A_i,1\leq i \leq M-1}\ = & d(\A)\Prob_{\A_0=\A}\set{x(z)=i,
\forall z\in A_{i,1}, 1\leq
i \leq M-1}\\\\
 = & d(\A)\Prob_{\A_0=\A}\set{A_{i,1}\subseteq R_x^i,\ 1\leq i
\leq M-1},
\end{array}\end{equation}\normalsize
where $R_x^i:=\set{z\in\Z:\ x(z)=i}$. Note that given $x$ and $\A$,
then the left hand side can be written as the product
\begin{equation*}\hspace{-1cm}
\begin{array}{rl} \Prob_{\eta_0=x}\set{\eta_1(z)=i,\forall z\in A_i,1\leq i \leq
M-1} =&
\displaystyle\prod_{i=1}^{M-1}\Prob_{\eta_0=x}\set{\eta_1(z)=i,\forall
z\in A_i}\\\\
 = &\displaystyle \prod_{i=1}^{M-1}\displaystyle\prod_{z\in A_i}\Prob_{\eta_0=x}\set{\eta_1(z)=i},
\end{array}
\end{equation*}
as long as $A_1,\ldots,A_{M-1}$ are pairwise disjoint. This
factorization, due to the conditional independence of different
particles, leads us to consider dual chains where the sets
$A_{1,s},\ldots,A_{M-1,s}$ evolve independently as long as they are
disjoint; that is, given pairwise disjoint $A_1,\ldots,A_{M-1}$, we
require
\begin{equation*}\Prob_{\A_0=\A}\set{A_{i,1}=B_i:\ 1\le i\le
M-1}=\prod_{i=1}^{M-1}\Prob_{\A_0=\A}\set{A_{i,1}=B_i}\end{equation*}
for all pairwise disjoint $B_1,\ldots,B_{M-1}\in\Y$. Indeed, we will
consider processes in which the evolution of points belonging only
to one component is independent of the rest. More precisely, for
pairwise disjoint $A_{1,s},\ldots,A_{M-1,s}$, a singleton $\{z\}$
contained in $A_{m,s}$ is substituted by $B_m(z)$, a subset of
$\set{z-1,z,z+1}$, independently of the rest of singletons in
$A_{m,s}$, with probabilities $\pi_m^\cdot$ given by

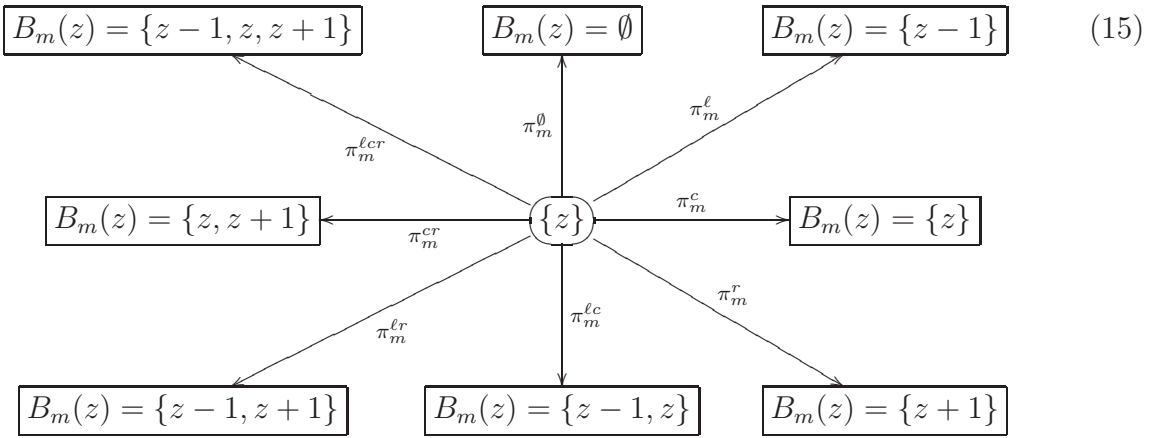
\begin{equation}\label{GraphDualTransitions}\hspace{-2cm}
\xymatrix{
      *+[F]{ B_m(z)=\set{z-1,z,z+1} }          &    *+[F]{  B_m(z)=\emptyset }     &      *+[F]{
      B_m(z)=\set{z-1} }
       \\\\
    *+[F]{  B_m(z)=\set{z,z+1} }           &   *+[F-:<10pt>]{   \set{z}}\ar[uu]^-{\pi_m^\emptyset}\ar[ruu]^-{\pi_m^\ell}
                 \ar[r]^-{\pi_m^c}\ar[rdd]^-{\pi_m^r}\ar[dd]^-{\pi_m^{\ell c}}
                 \ar[ldd]^-{\pi_m^{\ell r}}\ar[l]^-{\pi_m^{cr}}\ar[luu]^-{\pi_m^{\ell cr}}
                        &     *+[F]{   B_m(z)=\set{z}   }    \\\\
   *+[F]{      B_m(z)=\set{z-1,z+1} }       &  *+[F]{  B_m(z)=\set{z-1,z} }       &    *+[F]{
  B_m(z)= \set{z+1}}
}
\end{equation}

The set $A_{m,s+1}$ is the union of the resulting sets from the
updating of all $\{z\}\subseteq A_{m,s}$, that is
$A_{m,s+1}:=\bigcup_{z\in A_{m,s}}B_m(z)$. If the sets $A_i$ are not
pairwise disjoint, the point $\A=(A_1,\ldots,A_{M-1})$ is absorbing
for the dual process. Note that the evolution of this process can be
seen as a multi-type coalescing branching process. Moreover, we will
take the function $d$ with the form

\begin{equation}\label{d_votermodel1}
d(\A)=\prod_{i=1}^{M-1}d_i^{\abs{A_i}},
\end{equation}
with $d_i\in[0,+\infty)$.

Now, let $x\in\X$ and $\A\in\Y$. If $A_1,\ldots,A_{M-1}$ are not
pairwise disjoint, then both sides of (\ref{MONVM}) are zero. If
$A_1,\ldots,A_{M-1}$ are pairwise disjoint, we have, due to the
evolution of the dual, and noting that $R_x^1,\ldots,R_x^{M-1}$ are
pairwise disjoint by definition,
\begin{equation*}
\begin{array}{rl}\hspace{-2cm}
d(\A)\Prob_{\A_0=\A}\set{A_{i,1}\subseteq R_x^i:\ 1\le i\le M-1}
=&\displaystyle \prod_{i=1}^{M-1}d_i^{\abs{A_i}}\Prob_{\A_0=\A}\set{A_{i,1}\subseteq R_x^i}\\\\
=&
\displaystyle\prod_{i=1}^{M-1}d_i^{\abs{A_i}}\displaystyle\prod_{\{z\}\subseteq
A_i}\Prob_{\A_0=\A}\set{B_i(z)\subseteq R_x^i},
\end{array}
\end{equation*}
where $B_i(z)$ is the one-step evolution of $\{z\}\subseteq A_i$ in
the dual process.

Therefore, equation (\ref{MONVM}) takes the form
\begin{equation*}\prod_{i=1}^{M-1}\prod_{z\in
A_i}\Prob_{\eta_0=x}\set{\eta_1(z)=i}=\prod_{i=1}^{M-1}\prod_{\{z\}\subseteq
A_i} d_i\Prob_{\A_0=\A}\set{B_i(z)\subseteq R_x^i}.\end{equation*}

Now, we look for conditions on the probabilities $\pi_m^\cdot$ such
that $\eta_t$ and $\A_t$ are dual with respect to $(H,d)$. We must
have
$\Prob_{\eta_0=x}\set{\eta_1(z)=m}=d_m\Prob_{\A_0=\A}\set{B_m(z)\subseteq
R_x^m}$ for all $x\in\X$, $\A=(A_1,\ldots,A_{M-1})$ pairwise
disjoint, $z\in\Z$ and $m=1,\ldots,M-1$. These equations take the
following form
\begin{equation}\label{dualeq1}
\begin{array}{lcl}
p_{ijk,m} & = & d_m\pi_m^\emptyset \\
p_{mjk,m} & = & d_m\CC{\pi_m^\emptyset + \pi_m^\ell}\\
p_{imk,m} & = & d_m\CC{\pi_m^\emptyset + \pi_m^c}\\
p_{ijm,m} & = & d_m\CC{\pi_m^\emptyset + \pi_m^r}\\
p_{mmk,m} & = & d_m\CC{\pi_m^\emptyset + \pi_m^\ell + \pi_m^c +
\pi_m^{\ell c}}\\
p_{mjm,m} & = & d_m\CC{\pi_m^\emptyset + \pi_m^\ell + \pi_m^r +
\pi_m^{\ell r}}\\
p_{imm,m} & = & d_m\CC{\pi_m^\emptyset + \pi_m^c + \pi_m^r +
\pi_m^{c r}}\\
p_{mmm,m} & = & d_m\CC{\pi_m^\emptyset + \pi_m^\ell + \pi_m^c +
\pi_m^r + \pi_m^{\ell c} + \pi_m^{\ell r} + \pi_m^{c r}  +
\pi_m^{\ell c r}}
\end{array}
 ,\end{equation}
with  $i,j,k\in W\setminus\{m\}$, together with
\begin{equation}\label{pi}
\pi_m^\cdot\ge0,\qquad\qquad\pi_m^\emptyset + \pi_m^\ell + \pi_m^c +
\pi_m^r + \pi_m^{\ell c} + \pi_m^{\ell r} + \pi_m^{c r} +
\pi_m^{\ell c r}=1.
\end{equation}

Note that the last line of (\ref{dualeq1}) is actually

\begin{equation*}
p_{mmm,m}=d_m.
\end{equation*}

The existence of solutions for equations (\ref{dualeq1}) and
(\ref{pi}) requires extra assumptions on the process $\eta_t$. In
fact, such equations may only have solutions if for any $m\in W
\setminus\{M\}$ each of the parameters
$p_{ijk,m},p_{mjk,m},p_{imk,m},\ldots,p_{imm,m}$ does not depend on
$i,j,k\neq m$. That is to say that the probability of a site
assuming a state $m$ at time $s+1$ is a function of the positions of
the state $m$ in its neighbourhood at time $s$. Moreover, when
(\ref{dualeq1}) has a solution, we can note that $p_{mmm,m}=0$
implies $p_{ijk,m}=0$ for all $i,j,k\in W$, so the state $m$ could
not be taken by $\eta_t$. Hence, it is reasonable to assume
$p_{mmm,m}>0$ for all $m\in W\setminus\{M\}$. Thus, $0<d_m\le1$ for
all $m=1,\ldots,M-1$.

Therefore,  for any $m\in W\setminus\{M\}$ and $i,j,k\neq m$  we can
denote $p_m:=p_{ijk,m}$, which can be interpreted as the probability
that $\eta_1(z)$ spontaneously assumes the state $m$ (i.e. when $m$
is not present in the neighbourhood of the point). Thus, $p_{mmm,m}$
is equal to the probability that $\eta_1(z)$ does not spontaneously
assume any state $k\neq m$, and for $m=1,\ldots, M-1$,
 we have
\begin{equation*} d_m=1-p_{mmm,M}-p+p_m
\end{equation*}
where $p:=\sum_{n=1}^{M-1}p_n$. Hence, (\ref{d_votermodel1}) becomes
\begin{equation}\label{d_votermodel1b}
d(\A)=\prod_{i=1}^{M-1}(1-p_{iii,M}-p+p_i)^{\abs{A_i}}.
\end{equation}

The equations (\ref{dualeq1}) and (\ref{pi}) have solutions given
by: \footnotesize
\begin{equation}\label{solutioneq1}\hspace{-2cm}
\begin{array}{lcl}
\pi_m^\emptyset & = & \displaystyle\frac{1}{d_m}p_m\\
\pi_m^\ell  & = & \displaystyle\frac{1}{d_m}( p_{mkk,m}-p_m)\\
\pi_m^c  & = &  \displaystyle\frac{1}{d_m}(p_{kmk,m}-p_m)\\
\pi_m^r  & = & \displaystyle\frac{1}{d_m}(p_{kkm,m}-p_m)\\
\pi_m^{\ell c}  & = & \displaystyle
\frac{1}{d_m}(p_{mmk,m}+p_m-p_{mkk,m}-p_{kmk,m})\\
\pi_m^{\ell r}  & = & \displaystyle
\frac{1}{d_m}(p_{mkm,m}+p_m-p_{mkk,m}-p_{kkm,m})\\
\pi_m^{cr}  & = &  \displaystyle\frac{1}{d_m}(p_{kmm,m}+p_m-p_{kmk,m}-p_{kkm,m})\\
\pi_m^{\ell cr}  & = & \displaystyle\frac{1}{d_m}(p_{mmm,m}+p_{mkk,m}+p_{kmk,m}+p_{kkm,m}-p_{kmm,m}-p_{mkm,m}-p_{mmk,m}-p_m)
\end{array}
\end{equation}\normalsize
with $k\neq m$, provided that
\begin{equation}\label{ineq1}\hspace{-2cm}
\begin{array}{l}
 p_{mkk,m}\geq p_m\\
 p_{kmk,m}\geq p_m\\
 p_{kkm,m}\geq p_m\\
 p_{mmk,m} \geq p_{mkk,m}+p_{kmk,m}-p_m\\
 p_{mkm,m} \geq p_{mkk,m}+p_{kkm,m}-p_m\\
 p_{kmm,m} \geq p_{kmk,m}+p_{kkm,m}-p_m\\
 p_{mmm,m}+p_{mkk,m}+p_{kmk,m}+p_{kkm,m}\geq p_{kmm,m}+p_{mkm,m}+p_{mmk,m}+p_m
\end{array}
\end{equation}
for all $m=1,\ldots,M-1$, $k\in W\setminus\{m\}$. The inequalities
(\ref{ineq1}) impose lower and upper bounds on the velocity for the
growth of the probability that a site assumes the opinion $m$ in the
time $s+1$ because of the occurrence of opinions $m$ in its
neighbourhood in time $s$.\\

Thus, we have proved:

{\theo\label{votertheorem1} Let $\eta_t$ be a nearest-neighbour PCA
with state space $\X$, whose transition probabilities are such that
for any $m=1,\ldots,M-1$ each of the parameters
$p_{ijk,m},p_{mjk,m},p_{imk,m},\ldots,p_{imm,m}$ does not depend on
$i,j,k\neq m$ and satisfy inequalities (\ref{ineq1}). Then there
exist a Markov chain $\A_t$ with state space $\Y$,
$H:\X\times\Y\to\R$, and $d:\Y\to (0,1]$, such that $\eta_t$ and
 $\A_t$ are dual to one another with respect to
$(H,d)$. Furthermore, if $p_{mmm,m}=1$, for all $m=1,\ldots, M-1$,
then $\eta_t$ and  $\A_t$ are dual to one another with respect to
$H$.}

\cqd

The following result follows from the application of Theorem
\ref{theo(H,d)} and Theorem
\ref{votertheorem1}.\\

{\theo\label{corvoter1} Let $\eta_t$ be a nearest-neighbour PCA
whose transition probabilities are such that for any
$m=1,\ldots,M-1$ each of the parameters
$p_{ijk,m},p_{mjk,m},p_{imk,m},\ldots,p_{imm,m}$ does not depend on
$i,j,k\neq m$ and satisfy inequalities (\ref{ineq1}), and let
$\A_t=(A_{i,t})_{1\leq i\leq M-1}$ be the process on $\Y$ defined by
the transition probabilities (\ref{GraphDualTransitions}) and
(\ref{solutioneq1}). Then, the process is ergodic if and only if any
of the following conditions holds:
\begin{enumerate}

\item $p_{iii,M}>0$ for all $i\in W\setminus\{M\}$, that is $\eta_t$ admits spontaneous changes from any pure states $i\neq M$ to the state
$M$;

\item $\exists m,n\in W\setminus\{M\}$, $m\neq n$, with $p_m,p_n>0$, that is
$\eta_t$ admits spontaneous changes to at least two distinct states;

\item $\exists m\in W\setminus\{M\}$, with
$p_m>0$ and $p_{mmm,M}>0$;

\item $\exists m\in W\setminus\{M\}$, with $p_m>0$, and such that
the process $\A_t$ starting at
$\A=(\emptyset,\ldots,A_m,\ldots,\emptyset)$ has probability $1$ of
extinction for all $A_m\in Y$.
\end{enumerate}

Furthermore, when any of the above conditions holds, the unique
equilibrium measure for $\eta_t$ is defined for any
 $\A=(A_1,\ldots,A_{M-1})\in\Y$, by
\begin{equation*}
 \hat{\nu}(\A)=\Prob_{\tilde{\A}_{0}=\A}
\set{\tilde{\A}_{\tau}=  (\emptyset,\ldots,\emptyset)},
\end{equation*}
where $\tilde{\A}_t$ is the process defined on $\Y\cup\set{\wp}$
with transition probabilities given in (\ref{tilde xi_t}), and
$\tau$ is the hitting time of $
\{(\emptyset,\ldots,\emptyset)\}\cup\set{\wp}$ for $\tilde{\A}_t$.}

\proof We first check the sufficiency of the conditions for
ergodicity.

With the previous construction, $\eta_t$ and $\A_t$ are dual to one
another with respect to $(H,d)$. Also, by Lemma 3.1 of \cite{ls},
the set of linear combinations of functions $H(\cdot,\A)$ is dense
in $C(\X)$.

The set of absorbing states of $\A_t$  contains
$(\emptyset,\ldots,\emptyset)$ together with all
$\A=(A_1,\ldots,A_{M-1})$ which are not pairwise disjoint. Moreover,
let us denote ${\cal A}$ the set of values $m\in\{1,\ldots,M-1\}$,
if any, such that
\begin{equation*}\hspace{-2cm} p_m=p_{mkk,m}=p_{kkm,m}=p_{mkm,m}=0,\qquad
p_{kmk,m}=p_{mmk,m}=p_{kmm,m}=p_{mmm,m}.\end{equation*}

Then, if $m\in{\cal A}$, we have $\pi_m^c=1$ so the points in $A_m$
do not change in time. Thus, the set $\Theta$ of absorbing points of
$\A_t$ also includes the points $\A=(A_1,\ldots,A_{M-1})$ such that
$A_i=\emptyset$ for all $i\not\in{\cal A}$.

Let us see now that under any of the conditions {\em i}, {\em ii},
{\em iii} or {\em iv}, it holds $\{\theta\in\Theta:d(\theta)=1\}=
\{(\emptyset,\ldots,\emptyset)\}$. Obviously,
$d((\emptyset,\ldots,\emptyset))=1$. Now, for $(A_1,\ldots,A_{M-1})$
not pairwise disjoint with $A_i\cap A_j\ne\emptyset$, we have, under
any of {\em i}, {\em ii}, {\em iii} or {\em iv}, that at least one
of $d_i=1-p_{iii,M}-p+p_i$ or $d_j=1-p_{jjj,M}-p+p_j$ must be
smaller than 1, so from (\ref{d_votermodel1b}) we get $d(\A)<1$.
Last, let $i\in{\cal A}$; since $p_i=0$ we have that
$d_i=1-p_{iii,M}-p$; note that under {\em i} we have $p_{iii,M}>0$,
and under any of {\em ii}, {\em iii} or {\em iv}, there must be some
$m\in W\setminus(\{M\}\cup\mathcal{A})$ with $p_m>0$; therefore
under any of {\em i}, {\em ii}, {\em iii} or {\em iv}, we get
$d_i=1-p_{iii,M}-p\leq 1-p_{iii,M}-p_m<1$ and $d(\A)<1$. Note that
in particular, $\{\theta\in\Theta:d(\theta)=1\}=
\{(\emptyset,\ldots,\emptyset)\}$ implies that the condition {\em
iii} of Theorem \ref{theo(H,d)} holds with $c(\theta)=1$.

Suppose {\em i} holds, therefore it follows $\displaystyle d(\A)\leq
\max_{k\in W\setminus\{M\}}\{1-p_{kkk,M}\}<1$ for any non-empty set
$\A\in\Y$. If {\em ii} holds, then $\displaystyle d(\A)\leq
\max_{k=m,n}\{1-p_k\}<1$ for any non-empty set $\A\in\Y$. If {\em
iii} holds, then for any non-empty set $\A\in\Y$ we have $d(\A)\leq
\max\{1-p_m,\ 1-p_{mmm,M}\}<1$. Hence, in all the cases, we apply
Theorem \ref{theo(H,d)} to conclude.

Now, suppose {\em iv} holds. We can consider that for any $k\in
W\setminus\{M\}$, $k\neq m$, we have $p_k=0$, and also
$p_{mmm,M}=0$, since if the contrary holds, we could conclude the
ergodicity by {\em ii} or {\em iii}. Hence, for any $k\neq m$, it
follows that $\pi_k^\emptyset=0$, while $\pi_m^\emptyset>0$ (see
(\ref{solutioneq1})). Notice that $\Y=\Y_1\cup\Y_2$, where
$\Y_1:=\set{\B\in\Y:\ \exists k\neq m, B_k\neq\emptyset}$ and
$\Y_2:=\set{\B\in\Y:\ \forall k\neq m, B_k=\emptyset}$. For all
$\A\in\Y_1$ we have $d(\A)\le1-p_m<1$, and if $\A_0=\A$, then
$\A_s\in\Y_1$ for all $s<\tau$. Thus, the process $\tilde{\A}_t$
starting in $\A\in\Y_1$ has probability $1$ of being absorbed in
$\wp$. On the other hand, for all $\A\in\Y_2$ we have $d(\A)=1$ and
if $\tilde{\A}_t$ starts in $\A\in\Y_2$, then we have
$\tilde{\A}_s\in\Y_2$ for all $s<\tau$, and it coincides with the
process $\A_t$ starting in $\A\in\Y_2$ which has probability $1$ of
being absorbed in $(\emptyset,\ldots,\emptyset)$. Therefore, in both
cases, $\Prob_{\tilde\A_0=\A}\set{\tau<\infty}=1$. Thus, from the
proof of Theorem \ref{theo(H,d)} we conclude the ergodicity of
$\eta_t$ (see Remark \ref{remtau}).

Furthermore, in all the cases we get the expression of the
equilibrium measure for $\eta_t$ from Theorem \ref{theo(H,d)}.

For necessity we only have to show that the process is not ergodic
in the two following cases: (a) $p_i=0$ for all $i=1,\ldots,M-1$ and
there exists $j\neq M$ such that $p_{jjj,M}=0$; or (b) $ \exists!
m\in\set{1,\ldots,M-1}$ with $p_m>0$, $p_{mmm,M}=0$, and there is
some $A_m\in Y$ such that the process $\A_t$ starting in
$\A_0=(\emptyset,\ldots,A_m,\ldots,\emptyset)$ has probability
strictly less than $1$ of extinction. Case (a) is immediate since
the probability measures $\delta_M$ and $\delta_j$, concentrated on
pure state configurations $x'\in\X$ and $x''\in\X$, $x'(z)=M$ and
$x''(z)=j$ for all $z\in\Z$, respectively, are both equilibrium
measures for the process.

For case (b), note first that the measure $\delta_m$ is an
equilibrium measure for the process. Now let $x\in\X$ with $x(z)\neq
m$ for all $z\in\Z$ and $\A=(\emptyset,\ldots,A_m,\ldots,\emptyset)$
such that $\A_t$ has a probability less than $1$ of extinction. Note
that the set $\Y(m):=\set{\B\in\Y:B_i=\emptyset,\forall i\neq m}$ is
closed for the evolution of the dual process (that is, $\A_t\in
\Y(m)$ for all $t\ge 0$), and $d(\B)=1$ for all $\B\in\Y(m)$. Then,
for the process $\A_t$ enclosed in $\Y(m)$, the duality relation
(\ref{duald}) takes the form of (\ref{dual}), and for the particular
value of $x$ and $\A$ we have
\begin{equation*}
\begin{array}{rl}\Prob_{\eta_0=x}\set{\eta_s(z)=m,\ \forall z\in
A_m}=&\Prob_{\A_0=\A}\set{A_{m,s}\subseteq R_x^m}\\\\
=& \Prob_{\A_{0}=\A}\set{A_{m,s}=\emptyset}\\\\
\le & \Prob_{\A_{0}=\A}\set{\tau<\infty}<1.\end{array}
\end{equation*}

Then, $\limsup_{ s\to\infty}\Prob_{\eta_0=x}\set{\eta_s(z)=m,\
\forall z\in A_m}<1$ so the process starting with a measure with
support in $\set{x\in\X:\ x(z)\neq m, \forall z\in\Z}$ is not
converging to $\delta_m$ and the process is not ergodic. \cqd

{\rem Theorems \ref{votertheorem1} and \ref{corvoter1} can be
extended to a wide class of PCA with finite range interaction on the
integer lattice $\Z^d$ (\ref{PCA}), which are such that
$\Prob\set{\eta_{s+1}(z)=w|\ \eta_s(z+n)=w_n,\ \forall n\in N }$
does not depend on $w_n\neq w$. In fact, for any $d$ and
$N\subset\Z^d$ we can define a dual process $\xi_t$ on $\Y=Y^{M-1}$,
where $Y$ is the collection of finite subsets of $\Z^d$, and a point
$z$ belonging to $A_{m,s}$ is substituted by a subset of $\set{z+n:\
n\in N}$. Thus, in an analogous way to that in (\ref{ineq1}) it is
possible to find inequalities that delimit the class of PCA
for which we can extend our results.}\\

The PCA model presented in this subsection recover many
one-dimensional models that have been referred in the literature. In
fact, our model recovers $M$-state linear and non-linear PCA (which
were studied in the $2$-state version by \cite{coxdurrett} and
\cite{liggett2}), presenting new results for them.\\

{\ex\label{monvm} {\bf (Multi-opinion noisy voter model)} Suppose
$\eta_t$ is defined as follows: For each $m\in W=\{1,\ldots,M\}$,
take $q_m\in [0,1]$, and define $q=\sum_{m\in W}q_m$; set
$\alpha,\beta,\gamma\ge0$, such that $\alpha+\beta+\gamma=1-q$;
finally suppose $\eta_t$ has transition probabilities given by
$$p_{ijk,m}:=\alpha\mathbf{1}_{\{i\}}(m)+\beta\mathbf{1}_{\{j\}}(m)+\gamma\mathbf{1}_{\{k\}}(m)+q_m.$$

Note that for any $m\in W$ and $i,j,k\in W\setminus\{m\}$ we have
$q_m=p_{ijk,m}=:p_m$ (which here is also defined for $m=M$). Using
the language of voter models, $q_m$ is the probability of a voter
assuming the opinion $m$, independently of her own opinion, as well
as the opinions of her neighbours in the previous step, while
$\alpha$, $\beta$ and $\gamma$ are the weights given by each voter
to her own opinion and her neighbours' opinions.}

\cor\label{cormonvm} Let $\eta_t$ be a multi-opinion noisy voter
model. Then there exists a coalescing random walk $\A_t$ on the
family of all finite subsets of $\Z$ such that $\eta_t$ and $\A_t$
are dual to one another with respect to $(H,d)$ defined in
(\ref{H1}) and (\ref{d_votermodel1b}). Moreover, $\eta_t$ is ergodic
if and only if there exists $m\in\{1,\ldots,M\}$ such that $p_m>0$.

\proof It is easy to check that for any $m=1,\ldots,M-1$ the
parameters $p_{ijk,m},p_{mjk,m}$, $p_{imk,m},\ldots,p_{imm,m}$ do
not depend on $i,j,k\neq m$ and satisfy inequalities (\ref{ineq1}).
Hence, we can apply Theorem \ref{votertheorem1} to deduce the
existence of a dual with respect to $(H,d)$, where $H$ is given by
(\ref{H1}) and $d$ is given by (\ref{d_votermodel1b}).

For the second part of the result, if $p_m=0$ for all $m$, then any
probability measure concentrated on a pure state configuration
$x\in\X$, $x(z)=m$ for all $m\in\Z$ and any $m\in W$, is an
equilibrium measure. For sufficiency, note that (\ref{ineq1}) holds
regardless of the order of the opinions, and thus we can set $M$ to
have $p_M>0$ and ergodicity follows from Theorem \ref{corvoter1}
{\em i}.

\cqd

{\ex\label{dkm} {\bf (Domany-Kinzel model)} If we consider the case
$W=\set{0,1}$ and such that the probability of a site $z$ assuming
$1$ at time $s+1$ is a function of the number of $1$s in the
neighbourhood $\{z-1,z+1\}$ at time $s$, then we have the
Domany-Kinzel model (\ref{DK}) with parameters

$$\begin{array}{lcl}
a_0 & := & p_{000,1}=p_{010,1};\\\\
a_1 & := & p_{100,1}=p_{110,1}=p_{011,1}=p_{001,1};\\\\
a_2 & := & p_{101,1}=p_{111,1}.
\end{array}$$

We recall that the Domany-Kinzel model was introduced in
\cite{domanykinzel} and has been extensively studied (see for
instance \cite{konno2}, \cite{atmanmoreira}, \cite{kemper} and
\cite{rss}) due to its useful applicability in percolation theory
and phase-transition theory among others.\\}

{\cor\label{cordkm} Let $\eta_t$ be the Domany-Kinzel model
(\ref{DK}) with parameters $0\leq a_0\leq a_1\leq a_2\leq 1$. Then
there exists a coalescing random walk $\A_t$ on the family of all
finite subsets of $\Z$, such that $\eta_t$ and $\A_t$ are dual to
one another with respect to the pair $(H,d)$ defined by (\ref{H1})
and (\ref{d_votermodel1b}). Furthermore, if one of the following
conditions holds

\begin{enumerate}

\item $a_0+a_2\ge2a_1$ and $a_2<1$;

\item $a_0+a_2<2a_1$ and $a_0>0$;

\item $a_0>0$ and $a_2<1$;

\item $a_0=0$, $a_1< 1/2$ and $a_2<1$;

\item $a_0>0$, $a_1> 1/2$ and $a_2=1$;

\end{enumerate}

then $\eta_t$ is ergodic.}

\proof If $a_1\leq\frac{1}{2}(a_0+a_2)$, then $\eta_t$ satisfies all
hypotheses of Theorem \ref{votertheorem1} (identifying the state 0
here with the the state $M$ in Theorem \ref{votertheorem1}). Hence,
we have that $\A_t$ is defined by the transition probabilities
$\pi^\emptyset:={a_0}/{a_2}$, $\pi^\ell=\pi^r:=(a_1-a_0)/{a_2}$, and
$\pi^{\ell r}:=(a_2-2a_1+a_0)/{a_2}$.

On the other hand, if $a_1> \frac{1}{2}(a_0+a_2)$, then we can
consider the description of the process by the probabilities of the
occurrence of $0$, that is the Domany-Kinzel model with parameters
\begin{equation*}a_i^*:=1-a_{2-i},\qquad i=0,1,2.\end{equation*}

It is easy to see that $0\leq a_0^*\leq a_1^*\leq a_2^*\leq 1$, and
$a_1^*< \frac{1}{2}(a_0^*+a_2^*)$. Therefore, we can again apply
Theorem \ref{votertheorem1} (now the state $1$ plays the role of the
state $M$ in Theorem \ref{votertheorem1}) and find $\A_t^*$, which
is defined by the transition probabilities $\pi_*^\emptyset
:=a_0^*/a_2^*$, $\pi_*^\ell=\pi_*^r:=(a_1^*-a_0^*)/a_2^*$, and
$\pi_*^{\ell r}:=(a_2^*-2a_1^*+a_0^*)/a_2^*$.\\

Furthermore, we can obtain sufficient conditions for ergodicity:

\begin{enumerate}

\item Identifying the
state $M$ with 0 in Theorem \ref{votertheorem1} and
$p_{111,0}=1-a_2>0$, the result follows from Theorem \ref{corvoter1}
{\em i}.

\item We now identify $M$ with 1 so $p_{000,1}=1-a_2^*=a_0>0$ and the
conclusion follows from Theorem \ref{corvoter1} {\em i}.

\item We identify $M$ with 0 or 1 depending on $(a_2+a_0)/2$
being greater or smaller than $a_1$. In both cases, for $i\neq M$ we
have $p_{iii,M}>0$ and the result follows from Theorem
\ref{corvoter1} {\em i}.

\item Suppose $a_0=0$, $a_1<1/2$, and $a_2<1$. We have two possible
cases:

\begin{enumerate}

\item If $a_1\le a_2/2$, then we can use the dual process $\A_t$
and, since $p_{111,0}=1-a_2>0$ and the state $0$ plays the role of
the state $M$ in the Theorem \ref{votertheorem1}, we can deduce the
ergodicity of $\eta_t$ from   Theorem \ref{corvoter1}.{\em i}.

\item If $a_1> a_2/2$, then consider the dual process $\A_t^*$
(that is, with $M=  1$, $p_{000,1}=0$, $p_0=1-a_2>0$), starting from
any $\A_0\in\Y$.
 Let us show
that $\A_t^*$ has probability $1$ of extinction.

Consider the branching process $\chi_t$ on $\Z$ defined as follows:
Let $\chi_s$ be the number of individuals in the $s$th generation of
some population; assume each individual independently will give rise
to $0$ offsprings with probability $\pi_*^\emptyset$, $1$ offspring
with probability $(\pi_*^\ell+\pi_*^r)$, or $2$ offsprings with
probability $\pi_*^{\ell r}$; then $\chi_{s+1}$ is the number of
individuals in the next generation, that is the total number of
offsprings generated by the individuals in the generation $s$. Since
$a_1<1/2$, it follows that the expected number of offsprings of each
individual is $(\pi_*^l+\pi_*^r)+2\pi_*^{lr}=2a_1<1$, which implies
the branching process $\chi_t$ vanishes with probability $1$ (see
\cite{AN2}). Identifying $\chi_t$ starting with $\vert \A_0\vert$
individuals as an upper bound on the number of individuals of
$\A^*_t$, we have that the dual process $\A_t^*$ has probability 1
of being absorbed in ${\mathbf\emptyset}$ and we can use   Theorem
\ref{corvoter1}.{\em iv} to deduce the ergodicity of $\eta_t$.

\end{enumerate}

\item  If $a_0>0$, $a_1> 1/2$ and $a_2=1$, then $a_0^*=0$,
$a_1^*<1/2$ and $a_2^*<1$, and therefore we get the result reasoning
as in {\em iv}.

\end{enumerate}
\cqd

{\rem\label{d_relevance1} As we said in the Introduction, the
inclusion of the function $d$ in our definition of duality
(\ref{duald}) allows, in some cases, to analyze processes which do
not have a dual with respect to a given function $H$ in terms of the
classical duality equation (\ref{dual}). In fact, from the
construction of Theorem \ref{votertheorem1}, we can see that if we
consider the classical duality (\ref{dual}) with respect to the
function $H$ defined in (\ref{H1}) we would obtain a system of
equations similar to (\ref{dualeq1}), but without the factor $d_m$
multiplying the right side of the equations. Therefore, for the
existence of a solution of the system, we would need the extra
hypothesis $p_{mmm,m}=1$ for all $m=1,\dots,M-1$. In particular,
$M\ge 3$ would imply that every state $m=1,\ldots,M-1$ should have
$p_m=0$ and $p_{mmm,M}=0$; while if $M=2$, we could have $p_1>0$ but
not $p_{111,2}>0$. Note that the additional hypothesis of
$p_{mmm,m}=1$ for all $m\neq M$ would imply that if $M\geq 3$, then
$\eta_t$ is not ergodic; and if $M=2$, then $\eta_t$ would be
ergodic if and only if $p_1$ was large enough in order for condition
of Theorem \ref{corvoter1} {\em iii} to hold.

In Example \ref{monvm}, if we had used the classical duality
equation (\ref{dual}), we would be restricted to consider models
with every $q_m=0$ (if $M\geq 3$) or with only one $q_m\neq 0$ (if
$M=2$). On the other hand, in Example \ref{dkm}, and due to symmetry
of states `0' and `1' in the Domany-Kinzel model, we could have
either $a_0>0$ or $a_2>0$, but not both $a_0$ and $a_2$ positive.}

\subsection{Multi-state monotone models}

In the previous subsection we showed that a necessary condition for
a process to have a dual with respect to $(H,d)$ defined in
(\ref{H1}) and (\ref{d_votermodel1b}) is that for any
$m=1,\ldots,M-1$ each of the parameters $p_{...,m}$ depends only on
the presence (and position) of the value $m$ in the neighbourhood of
the point. However, there are models (for instance, when the set $W$
has an order) where not only the presence (and position) of the
value $m$ is important for the probability of changing to $m$ but
also the presence (and position) of values bigger or smaller than
$m$ is. In order to find a dual for these processes, we will use a
duality function $H$ different to (\ref{H1}) and allow more
transitions for the dual than we did in Subsection \ref{ssmonvm}; as
we will see, in this case, the evolution of any two subsets $A_i$
and $A_j$ has a greater interdependence than in the non-biased case.

For simplicity, we consider only the case of $\eta_t$ being a
nearest-neighbour PCA for which the probability of a fixed site
assuming some state at time $s+1$ does not depend on its own state
at time $s$. That is, we assume that
\begin{equation}\label{Condition1BiasedModel}
p_{ik,m}:=p_{ijk,m},\qquad \forall i,j,k,m\in W.
\end{equation}

Throughout this subsection, the state space of the dual process
will be $\Y=\set{\A\in Y^{M-1}:\ A_1\subseteq\cdots\subseteq
A_{M-1}}$ and the function $H:\X\times\Y\to\R$, defined by
\begin{equation}\label{biasedH} H(x,\A)=\left\{\begin{array}{lll} 1  &,&\ {if}\
x(z)\leq k,\ \forall z\in A_k,\ \forall k=1,\ldots, M-1\\
0 &,&\ {otherwise}\end{array}\right..\end{equation}

The class of dual processes we consider is composed of processes
$\A_t(=\xi_t)$ starting at $\A_0\in\Y$ whose transition
probabilities from $\A_s=(A_{i,s})_{1\leq i\leq M-1}$ to
$\A_{s+1}=(A_{i,s+1})_{1\leq i\leq M-1}$ are defined as follows. Let
$(\emptyset,\ldots,\emptyset)$ be an absorbing state for $\A_t$ and,
under the convention that $A_{0,s}=\emptyset$, assume that for any
$k,m,n\in\{1,\ldots,M-1\}$ each singleton $\{z\}\in A_{k,s}\setminus
A_{k-1,s}$ will generate an element independently of other $  z'\in
A_{M-1,s}$, $\B(z)=(B_i(z))_{1\leq i\leq M-1}\in\Y$ as follows:
\begin{itemize} \item with probability $\pi_k^\emptyset$
$$B_i(z)=\emptyset,\qquad 1\leq i\leq M-1;$$

\item with probability $\pi_{k,m}^\ell$
$$B_i(z)=\left\{\begin{array}{lcl}
\emptyset & {,\ if } & 1\leq i < m\\\\
\{z-1\} & {,\ if } & m\leq i\leq M-1;
\end{array}\right.$$

\item with probability $\pi_{k,m}^r$
$$B_i(z)=\left\{\begin{array}{lcl}
\emptyset & {,\ if } & 1\leq i < m\\\\
\{z+1\} & {,\ if } & m\leq i\leq M-1;
\end{array}\right.$$

\item with probability $\pi_{k,mn}$
$$B_i(z)=\left\{\begin{array}{lcl}
\emptyset & {,\ if } & 1\leq i < \min\{m,n\}\\\\
\{z-1\} & {,\ if } & m\le i<n\\\\
\{z+1\} & {,\ if } & n\le i<m\\\\
\{z-1,z+1\} & {,\ if } & \max\{m,n\}\leq i\leq M-1;
\end{array}\right.$$
\end{itemize}
with
$$\pi_k^\emptyset+ \sum_{m=1}^{M-1}\pi_{k,m}^\ell+
\sum_{n=1}^{M-1}\pi_{k,n}^r+
\sum_{m=1}^{M-1}\sum_{n=1}^{M-1}\pi_{k,mn}=1.$$

We define $\A_{s+1}=(A_{i,{s+1}})_{1\leq i\leq M-1}$, where
$$A_{i,s+1}=\bigcup_{z\in A_{M-1}} B_i(z).$$

Note that $\A_{s+1}$ also belongs to $\Y$.

Equation (\ref{duald}) is, in this case, for each $x\in \X$,
$\A\in\Y$,

\footnotesize
\begin{equation}\label{dualdbm}\hspace{-2cm}
\Prob_{\eta_0=x}\set{\eta_1(z)\le k: \forall z\in A_k, 1\leq k\leq
M-1}=d(\A)\Prob_{\A_0=\A}\set{A_{k,1}\subseteq \mathcal R_x^k:1\le
k\le M-1}, \end{equation} \normalsize where $\mathcal
R_x^k:=\set{z\in\Z:\ x(z)\leq k}$.

The left hand side of (\ref{dualdbm}) takes the form:
$$\displaystyle\prod_{k=1}^{M-1}\prod_{z\in A_k\setminus
A_{k-1}}\Prob_{\eta_0=x}\set{\eta_1(z)\le k}$$ and, assuming that
$d:\Y\to\R$ can be written as
$$d(\A)=\displaystyle\prod_{k=1}^{M-1} d_k^{\vert A_k\setminus
A_{k-1}\vert}$$ for some $d_1,\ldots,d_{M-1}\in[0,1]$ (the product
over the empty set is taken as 1), the right hand side of
(\ref{dualdbm}) is
$$\displaystyle\prod_{k=1}^{M-1}\prod_{\{z\}\subseteq A_k\setminus
A_{k-1}}d_k\Prob_{\A_0=\A}\set{B_m(z)\subseteq \mathcal R_x^m: 1\le
m\le M-1}.$$

Then, we must find solutions for the set of equations

\begin{equation}\label{ecuacionesbm}\Prob_{\eta_0=x}\set{\eta_1(z)\le k}=
d_k\Prob_{\A_0=\A}\set{B_m(z)\subseteq \mathcal R_x^m: 1\le m\le
M-1},\end{equation} for all $x\in \X$, $\A\in\Y$, $k=1,\ldots,M-1$
and $z\in A_k\setminus A_{k-1}$. Writing
$S_{i,j}^k=\sum_{m=1}^kp_{ij,m}$ and noting that, for $1\le i,j\le
M$, $x(z-1)=i$, $x(z+1)=j$, we have
\begin{equation*}\mathcal
R_x^m\supseteq\left\{\begin{array}{lll} \set{z-1} &,&\ {if}\
i\le m<j \\
\set{z+1} &,&\ {if}\
j\le m<i \\
\set{z-1,z+1} &,&\ {for}\ m\ge i,j
\end{array}\right..\end{equation*}

  The equations (\ref{ecuacionesbm}) are

\begin{equation}\label{dualbiased}\hspace{-2.5cm}
\begin{array}{lcl}
S_{M,M}^k & = & d_k\pi_k^\emptyset\\\\
S_{M-1,M}^k  & = & d_k\CC{\pi_k^\emptyset+ \pi_{k,M-1}^\ell}\\\\
S_{M,M-1}^k  & = & d_k\CC{\pi_k^\emptyset+ \pi_{k,M-1}^r}\\\\
S_{M-1,M-1}^k  & = & d_k\CC{\pi_k^\emptyset+ \pi_{k,M-1}^\ell+ \pi_{k,M-1}^r+ \pi_{k,M-1\ M-1}}\\\\
S_{M-2,M-1}^k  & = & d_k\CC{\pi_k^\emptyset+ \pi_{k,M-1}^\ell+
\pi_{k,M-2}^\ell
+ \pi_{k,M-1}^r+ \pi_{k,M-1\ M-1}+ \pi_{k,M-2\ M-1}}\\\\
\vdots &&\\\\
S_{i,j}^k  & = &\displaystyle d_k\CC{\pi_k^\emptyset+
\sum_{m=i}^{M-1}\pi_{k,m}^\ell+ \sum_{n=j}^{M-1}\pi_{k,n}^r+
\sum_{m=i}^{M-1}\sum_{n=j}^{M-1}\pi_{k,mn}}\\\\
\vdots &&\\\\
S_{1,1}^k  & = &\displaystyle d_k\CC{\pi_k^\emptyset+
\sum_{m=1}^{M-1}\pi_{k,m}^\ell+ \sum_{n=1}^{M-1}\pi_{k,n}^r+
\sum_{m=1}^{M-1}\sum_{n=1}^{M-1}\pi_{k,mn}}
\end{array}
\end{equation}
for $k=1,\ldots, M-1$. Due to the last line of the above equations,
we deduce that $d_k=S_{1,1}^k$ and then
\begin{equation}\label{d_votermodel2}d(\A)=\prod_{k=1}^{M-1}
\left(\sum_{m=1}^k p_{11,m}\right)^{\vert A_k\setminus
A_{k-1}\vert}.\end{equation} We assume $p_{11,1}=S_{1,1}^1>0$, for
otherwise, under (\ref{dualbiased}), $S_{i,j}^1=0$ for all
$i,j\in\{1,\ldots,M\}$ and opinion 1 cannot be taken by the process.
Therefore, $d_k=S_{1,1}^k>0$ for all $k\in \set{1,\ldots,M-1}$ and
the solution of (\ref{dualbiased}) is given by:
\begin{equation}\label{soldualbiased}
\begin{array}{lcl}
\pi_k^\emptyset & = & \displaystyle \frac{S_{M,M}^k}{d_k}\\\\
\pi_{k,m}^\ell & = & \displaystyle \frac{S_{m,M}^k-S_{m+1,M}^k}{d_k}\\\\
\pi_{k,m}^r & = & \displaystyle \frac{S_{M,m}^k-S_{M,m+1}^k}{d_k}\\\\
\pi_{k,mn} & = & \displaystyle
\frac{S_{m,n}^k-S_{m+1,n}^k-S_{m,n+1}^k+S_{m+1,n+1}^k}{d_k}
\end{array}
\end{equation}
since

\begin{equation}\label{solbm}\hspace{-1cm}
\begin{array}{lcl}
\displaystyle
d_k\sum_{m=i}^{M-1}\pi_{k,m}^\ell&=&\displaystyle\sum_{m=i}^{M-1}\left(S_{m,M}^k-S_{m+1,M}^k\right)
=S_{i,M}^k-S_{M,M}^k,\\\\
d_k\displaystyle\sum_{n=j}^{M-1}\pi_{k,n}^r&=&\displaystyle\sum_{n=j}^{M-1}\left(S_{M,n}^k-S_{M,n+1}^k\right)
=S_{M,j}^k-S_{M,M}^k,\\\\
d_k\displaystyle\sum_{m=i}^{M-1}\displaystyle\sum_{n=j}^{M-1}\pi_{k,mn}&=&
\displaystyle\sum_{m=i}^{M-1}\displaystyle\sum_{n=j}^{M-1}\left(S_{m,n}^k-S_{m+1,n}^k-S_{m,n+1}^k+S_{m+1,n+1}^k\right)\\\\
&=&\displaystyle\sum_{m=i}^{M-1}\left(S_{m,j}^k-S_{m,M}^k-S_{m+1,j}^k+S_{m+1,M}^k\right)\\\\
&=& S_{i,j}^k-S_{M,j}^k-S_{i,M}^k+S_{M,M}^k.
\end{array}
\end{equation}

In order to have the nonnegativity of (\ref{soldualbiased}), we need

\begin{equation}\label{cond1bm}S_{m,M}^k\ge S_{m+1,M}^k,\qquad\qquad S_{M,n}^k\ge
S_{M,n+1}^k,\end{equation}

\begin{equation}\label{cond2bm}
S_{m,n}^k-S_{m+1,n}^k-S_{m,n+1}^k+S_{m+1,n+1}^k\ge0,\end{equation}
for all $1\le m,n\le M-1$. Note that (\ref{cond1bm}) and
(\ref{cond2bm}) together imply $S_{i,j}^k\ge S_{m,n}^k$ for all
$1\le i\le m\le M$, $1\le j\le n\le M$ and $k=1,\ldots,M$, that is

\begin{equation*}\sum_{\ell=1}^kp_{ij,\ell}\ge\sum_{\ell=1}^kp_{mn,\ell}\end{equation*}
which means that the process $\eta_t$ is {\bf monotone} with respect
to the coordinatewise order in $\X$ induced by the total order
$1\le2\le\cdots\le M$ in $W$ (see \cite{LS00}). Therefore, condition
(\ref{cond1bm}) can be substituted by the monotonicity of the
process.

On the other hand, conditions (\ref{cond1bm}) and (\ref{cond2bm})
can be written equivalently as saying that $S_{i,j}^k=1-F^k(i,j)$
for each $k=1,\ldots M-1$, where $F^k$ is a distribution function on
$\R^2$.

We have proved the following result:

{\theo\label{votertheorem3} Let $\eta_t$ be a nearest-neighbour PCA
with state space $\X$ whose transition probabilities satisfy
(\ref{Condition1BiasedModel}). Then, there exist a dual process
$\A_t$ with state space $\Y$ such that $\eta_t$ and $\A_t$ are dual
to one another with respect to $(H,d)$ defined in (\ref{biasedH})
and (\ref{d_votermodel2}) if and only if $\eta_t$ is monotone and
its transition probabilities satisfy (\ref{cond2bm}). Furthermore,
if $p_{11,1}=1$, then $\eta_t$ and $\A_t$ are dual to one another
with respect to $H$.}

  \cqd

{\cor\label{corvoter3} A monotone nearest-neighbour PCA whose
transition probabilities satisfy (\ref{cond2bm}) is ergodic if
$p_{11,M}>0$.   In such case, the unique equilibrium measure for
$\eta_t$ is defined for any
 $\A=(A_1,\ldots,A_{M-1})\in\Y$, by
\begin{equation*}
 \hat{\nu}(\A)=\Prob_{\tilde{\A}_{0}=\A}
\set{\tilde{\A}_{\tau}=  (\emptyset,\ldots,\emptyset)},
\end{equation*}
where $\tilde{\A}_t$ is the process defined on $\Y\cup\set{\wp}$
with transition probabilities given in (\ref{soldualbiased}), and
$\tau$ is the hitting time of $
\{(\emptyset,\ldots,\emptyset)\}\cup\set{\wp}$ for $\tilde{\A}_t$.}

\proof

By Theorem \ref{votertheorem3} there exist a Markov chain $\A_t$
and a pair of functions $(H,d)$, such that  $\eta_t$ and $\A_t$
are dual to one another with respect to $(H,d)$.

We now apply Theorem \ref{theo(H,d)} to get the result. In fact,
by Lemma 3.2 of \cite{ls} the set of linear combinations of the
functions $H(\cdot,\A)$ with $\A\in\Y$ is dense in $C(\X)$.
  Also, $(\emptyset,\ldots,\emptyset)$ is the unique
absorbing state $\theta$ of $\A_t$ such that $d(\theta)=1$, and
$H(\cdot,(\emptyset,\ldots,\emptyset))\equiv1$. Moreover, for any
  $\A\neq(\emptyset,\ldots,\emptyset)$,
$d(\A)=\prod_{k=1}^{M-1}(1-\sum_{j=k+1}^M
p_{11,j})^{\abs{A_k\setminus A_{k-1}}}\le1-p_{11,M}$. Thus, if
$p_{11,m}>0$, hypotheses {\em ii} and {\em iii} (with
$c({\mathbf\emptyset})=1$) of Theorem \ref{theo(H,d)} hold and the
result is proved.
 \cqd

{\rem Note that as in the case of the multi-opinion voter models, we
can extend Theorem \ref{votertheorem3} and Corollary \ref{corvoter3}
to a wide class of multi-state monotone biased models with state
space $W^{\Z^d}$ and finite
range interaction.}\\

We notice that the processes $\eta_t$ that verify the conditions of
Theorem \ref{votertheorem3} include many PCA models, such as
multi-opinion noisy biased voter models and multi-species
competition models.

{\ex\label{mscm} {\bf (Multi-species competition model)} Suppose $M$
distinct species live in the space $\Z$. Each site of the space
contains only one individual. The spread mechanism for all species
is the same: each individual produces one offspring (independently
of other individuals) that will leave to colonize a neighbouring
site, killing the old particle that occupied the site and competing
against the offspring of others individuals. Moreover, we include
the possibility of a spontaneous appearance of an individual of any
type (for instance, coming from outside).

We consider such a model described as follows: Suppose $\eta_t$ is a
nearest-neighbour process with $M$ states such that the  probability
$p_{ij,m}$ that a site assumes the state $m$ in the next step when
its left neighbour is in the state $i$ and its right neighbour is in
the state $j$ is given for any $j,k\neq m$ and $i<m$ by
\begin{equation*}\begin{array}{l}
 p_{jk,m}= p_m \\\\
 p_{im,m}= p_m + (1-p)\alpha_m\\\\
 p_{im,i}= p_i + (1-p)(1-\alpha_m)\\\\
 p_{mi,m}= p_m + (1-p)\beta_m\\\\
 p_{mi,i}= p_i + (1-p)(1-\beta_m)\\\\
 p_{mm,m}= 1-p+p_m,
\end{array}
\end{equation*}
where $p=\sum_{m=1}^M p_m$, $\alpha_M\geq\alpha_{M-1}\ldots\geq
\alpha_2\geq 1/2$ and $\beta_M\geq\beta_{M-1}\ldots\geq
\beta_2\geq 1/2$.

The parameters $\alpha_m$ and $\beta_m$ can be interpreted as the
competitiveness of species $m$ with relation to species $i<m$, and
$p_m$ corresponds to probabilities that spontaneous mutations to $m$
occur. Note that the order in the $\alpha$s and $\beta$s imply that
species are ordered from weakest 1 to strongest $M$. Usually, we
have $\alpha_m=\beta_m$ for all $m$, which means there is not a
directional component in the competition.

 Competition models belong to the class of biased
voter models. The biased voter model was introduced in the $2$-state
version by Williams and Bjerknes \cite{willbjerk} as a model of skin
cancer, and has been applied to study several phenomena (for
instance, see \cite{jtp2007}, \cite{KordLey}, and \cite{Koulis}).

Since
$$S_{i,j}^k=\displaystyle\sum_{\ell=1}^kp_\ell+
(1-p)\left((1-\alpha_j){\bf 1}_{\{i\le k<j\}}+(1-\beta_i){\bf
1}_{\{j\le k<i\}}+{\bf 1}_{\{k\ge i,j\}}\right),$$ it is not hard to
check that $\eta_t$ is monotone and (\ref{cond2bm}) holds for it.
Therefore, by Corollary \ref{corvoter3}, if $\eta_t$ admits a
spontaneous change to the state $M$, then it is ergodic. }

{\ex\label{model_g} Consider the case of monotone processes for
which $S_{i,j}^k=g_k(i+j)$, where $g_k:\set{2,\ldots,2M}\to [0,1]$.
In such case, writing $l:=i+j$, condition (\ref{cond2bm}) takes the
form
$$\frac{1}{2}\C{g_k(l)+g_k(l+2)}\geq g_k(l+1),$$
which means that (\ref{cond2bm}) is satisfied if and only if the
functions $g_k$ are convex. In this case, we have that $\eta_t$ is
ergodic if $g_{M-1}(2)<1$. In fact, if $g_{M-1}(2)<1$, then
$p_{11,M}=1-S_{1,1}^{M-1}=1-g_{M-1}(2) >0$ and we can use Corollary
\ref{corvoter3} to deduce the ergodicity of the process. }

{\rem\label{d_relevance2} As in Subsection \ref{ssmonvm}, our
definition (\ref{duald}) of duality with respect to $(H,d)$ allows
us to include processes that do not have a dual in terms of the
classical definition (\ref{dual}). In fact, we could not obtain a
version of Theorem \ref{votertheorem3} for classical duality
(\ref{dual}) with respect to the function $H$ defined in
(\ref{biasedH}) unless we imposed another restrictive condition on
$\eta_t$ since we would obtain a system of equations similar to
(\ref{dualbiased}), but without the factor $d_k$ multiplying the
right side of each equation. Therefore, we should have $S_{1,1}^k=1$
for any $k=1\ldots,M-1$, which is equivalent to imposing the
condition $p_{11,1}=1$. In such case, we would have $p_{11,M}=0$ and
it would not be possible to obtain a sufficient condition for
ergodicity in the form of Corollary \ref{corvoter3}.

In particular, in Example \ref{mscm}, classical duality (\ref{dual})
with respect to $H$ would need  $p_m=0$ for any $m\neq 1$ (which
implies $p_M=0$, and therefore we could not deduce the ergodicity),
while in Example \ref{model_g} we should impose the condition
$g_1(2)=1$ (which implies $g_{M-1}(2)=1$) and therefore we could not
deduce the ergodicity either. }

\section{Concluding remarks}

In this work, we have presented a new duality equation (\ref{duald})
and shown its usefulness in the study of PCA. We believe that there
is still room for future research on this topic. In particular we
are considering the following issues:

\begin{enumerate}

\item Studying in more detail the equilibrium measures (proportion of
values, spatial correlations, etc.) when the processes are ergodic
in Examples \ref{monvm}, \ref{dkm} and \ref{mscm}.

\item In the case of processes having a dual with respect to
$(H,d)$ in Theorems \ref{votertheorem1} and \ref{votertheorem3} but
not being ergodic, studying their asymptotic behaviour using their
dual processes.

\item Studying other classes of functions $H$ giving conditions on
the transition probabilities for the existence of a dual with
respect to $(H,d)$.

\end{enumerate}

\ack

We thank the referees for their careful reading and useful comments.

This work was supported by projects MTM2004-01175 and MTM2007-63769
of MEC, and research group ``Modelos Estoc\'{a}sticos'' (DGA). M.
Sobottka was partially supported by Nucleus Millennium Information
and Randomness ICM P04-069-F and project DIUC 207.013.030-1.0.

\end{document}